# Gradient descent procedure for solving linear programming relaxations of combinatorial optimization problems in parallel mode on extra large scale

Alexey Antonov
aav.antonov@gmail.com

Linear programming (LP) relaxation is a standard technique for solving hard combinatorial optimization (CO) problems. Here we present a gradient descent algorithm which exploits the special structure of some LP relaxations induced by CO problems. The algorithm can be run in parallel mode and was implemented as CUDA C/C++ program to be executed on GPU. We exemplify efficiency of the algorithm by solving a fractional 2-matching (F2M) problem for instances of the traveling salesman problem (TSP) ranging in size from 100,000 up to 200,000 cities. Our results demonstrate that a fractional 2-matching problem with 100,000 nodes is solved by our algorithm on a modern GPU on a scale of a second while solving the same problem with simplex method would take more than an hour. The algorithm can be modified to solve more complicated LP relaxations derived from CO problems. The program code is freely available at https://github.com/aav-antonov/F2M .

## 1. Introduction

Linear programming (LP) relaxations are working horses for solving hard combinatorial optimization problems (Pulleyblank 1983). Commonly the LP relaxations are solved using standard LP technique (simplex method). Nevertheless being very efficient this technique is hard to parallelize and even with modern computational power is limited for solving problems on a modest scale. In practice, some LP relaxations with a special structure are solved by faster algorithms (max-flow algorithm or a heuristic) (Pulleyblank 1983, Schrijver 1998). However, adaptation of these algorithms for more complicated relaxations is not straightforward.

Here we propose to solve LPs that have a special structure based on general properties of linear programs. Namely, in case a LP has equality constraints one can modify the objective function (OF) of the LP by adding any linear combination of the equality constraints to the OF without actually affecting the optimal solution (the optimal solution would be the same for the original OF and for the modified one ). In some cases the solution of LP could directly follow from the OF (for example, for TSP, if edges of a given hamilton cycle have value 1 while all other edges have value 2). So instead of solving the original problem directly one can try to find a modification algorithm that would convert the original OF into a trivial one. In some cases for many LP relaxations of CO problems a gradient descent procedure can be developed to find a modified OF that leads to a trivial solution. In difference to the simplex method, the gradient descent is naturally parallelised and, therefore, could use modern GPUs to massively speed up computations.

Here we propose a gradient descent procedure (GDP) to solve a special class of LP relaxations induced by CO problems. Fractional 2-matching (F2M) is one of the common CO relaxation used for solving TSP and is an example of LP which can be solved efficiently by GDP. We have implemented GDP for F2M as CUDA C/C++ program to be executed on GPU. We exemplify efficiency of the algorithm by solving F2M on the scale of up 200,000 nodes. We demonstrate that an average F2M problem with 100,000 nodes is solved on a modern GPU on a scale of a second. We compare the performance of our program with a popular LP solver Soplex (Gleixner et al. 2018). The performance of Soplex to solve an average F2M problem with 100,000 nodes as a linear programming problem has been at least 1000 times slower.

## 2. Linear programming

Let us consider a general form of a LP problem:

$$\begin{aligned}&\textit{minimize } cx\\&\textit{subject to}: \qquad (LP1)\end{aligned}$$

$$a_i x = b_i \qquad i \in I \qquad (1)$$

$$a_i x \geq b_i \qquad i \in J \qquad (2)$$

$$x \geq 0 \qquad (3)$$

$$c, x \in R^n$$

Let $x_o \in R^n$ be the optimal solution of (LP1) and $Q$ is a subset of $J$ so that $a_q x_o = b_q$ for each $q \in Q \in J$. Consider a new objective function

$$z = c - \sum_{i \in I} \lambda_i a_i - \sum_{q \in Q} \lambda_q a_q \qquad (4)$$

where $\lambda$ is subjected to:

$$\lambda_q \geq 0 \; for\; q \in Q \;\; and\; \lambda_i \in R \; for\; i \in I \qquad (5)$$

From general properties of linear programs it is obvious that $x_o$ is also optimal for LP2 for any valid $\lambda$ subjected to (5):

$$\textit{minimize } zx = cx - \sum_{i \in I} \lambda_i (a_i x) - \sum_{q \in Q \in J} \lambda_q (a_q x)$$

$$\begin{aligned}&\textit{subject to}: \qquad (LP2)\\&(1), (2), (3),\\&c, x \in R^n\end{aligned}$$

The LP relaxations of CO problems commonly have additional binary constraint ( $0 \le x \le 1$ ). Let us denote as LP3 a class of LP problems with the structure presented below:

$$
\begin{aligned}
& \textit{minimize } cx \\
& \textit{subject to}: \qquad\qquad\qquad\qquad (LP3) \\
& \quad (1),\ (2), \\
& 0 \le x \le 1 \qquad\qquad\qquad\qquad (6) \\
& c, x \in R^n
\end{aligned}
$$

Assume $x_o$ is the optimal solution for LP3 and, therefore (assuming absence of degeneration), $x_o$ is a vertex of the polytope defined by (1), (2) and (6). From Karush–Kuhn–Tucker (KKT) optimality conditions it follows that objective function must be decomposed as linear combination of active constraints (1), (2) and (6) :

$$
c = \sum_{i \in I} \lambda_i a_i + \sum_{q \in Q} \lambda_q a_q + \sum_{k \in E^+} \lambda_+^k e_+^k + \sum_{k \in E^-} \lambda_-^k e_-^k \qquad \text{(KKT)}
$$

$$
\lambda_q > 0,\ \lambda_+^k > 0,\ \lambda_-^k > 0
$$

Vectors $e_+^k$ and vectors $e_-^k$ are unit vectors corresponding to active constraints (6) : $x_o^k = 1$ and $x_o^k =$ 0. Each vector $e_+^k = (0,\ ..,\ 1,\ ..0)$ has value 1 in position k and 0 otherwise while each vector $e_-^k = (0,\ ..,\ -1,\ ..0)$ has value -1 in position k and 0 otherwise. Therefore, vector

$$
v = c - \sum_{i \in I} \lambda_i a_i - \sum_{q \in Q} \lambda_q a_q
$$

due to KKT conditions must be decomposed as

$$
v = \sum_{k \in E^+} \lambda_+^k e_+^k + \sum_{k \in E^-} \lambda_-^k e_-^k
$$

and therefore we have :

$$
\begin{aligned}
& v^k < 0 \textit{ if } x_o^k = 1 \\
& v^k = 0 \textit{ if } (0 < x_o^k < 1) \qquad\qquad\qquad (7) \\
& v^k > 0 \textit{ if } x_o^k = 0
\end{aligned}
$$

Substitution of $v$ as objective function in place of $c$ in LP3 does not obviously change optimal solution $x_o$ and converts LP3 into a trivial problem (the optimal solution $x_o$ follows from $v$, see (7)).

## 3. Gradient descent

In many LP relaxations in addition to constraints (6) coefficients ( $a_i$ ) are also binary (1 or 0) while right sides ( $b_i$ ) are integral:

$$a_i \text{ is binary}, \; b_i \text{ is integral} \qquad (8)$$

Let us consider a constraint $a_i$ ( $i \in I \; or \; i \in Q$ ) and let $K_i$ to be the set of indexes so that $a_i^k = 1$ in case $k \in K_i$ and $a_i^k = 0$ otherwise ( $k \notin K_i$ ). Consider values $v^k$ $(k \in K_i)$ that correspond to optimal values $\lambda = (\lambda_i, \lambda_q)$. In this case the number of $v^k$ $(k \in K_i)$ with strictly negative values must be less or equal to $b_i$ ( $b_i$ is integral ). In other case constraint $i$ is violated as for each $v^k < 0$ we have $x_o^k = 1$ (see (7)). On the other hand, the number of $v^k$ with strictly positive values must be less or equal to $|K_i| - b_i$ by the same logic. This leads us to simple gradient formula for $\lambda_i$. Let us assume that we order $v^k$ $(k \in K_i)$ values in ascending order and $(v^k)_{b_i}$ , $(v^k)_{b_i+1}$ are elements in positions $b_i$ and $b_i + 1$ respectively. So to ensure (KKT for LP3) conditions the values of ( $(v^k)_{b_i}$ , $(v^k)_{b_i+1}$ ) must either satisfy to

$$(v^k)_{b_i} < 0 \text{ and } (v^k)_{b_i+1} > 0$$

or

$$(v^k)_{b_i} = (v^k)_{b_i+1} = 0$$

This lead to gradient for $\lambda_i$ to be:

$$grad(\lambda_i) = ((v^k)_{b_i} - (v^k)_{b_i+1})$$

Obviously, for optimal values $\lambda = (\lambda_i, \lambda_q)$ we have $grad(\lambda_i) = 0$ for all $i \in I \; and \; i \in Q$. However, the condition $grad(\lambda_i) = 0$ for all $i \in I \; and \; i \in Q$ does not guarantee the optimal solution to LP3. Our gradient descent procedure does not control exact value of constraints and may converge (in degeneration cases) to a different optimum with wrong value for some constraints: $a_i x = b_i \pm 1, 2, ..$ On the other hand, the computation of gradient does not require explicit value of $x$ thus making it very fast. There are different ways to implement control for violated constraints depending on the specificity of the problem in hand.

## 4. Fractional 2-matching

Fractional 2-matching (F2M) is a LP relaxation of the 2-matching problem which, in turn, is a relaxation of traveling salesman problem (TSP). TSP is formulated on an undirected graph $G = (V, E)$ where V is the set of graph nodes and E is a set of graph edges (Boyd and Carr 1999). For each edge $e \in E$ the weights $c_e$ are given and the aim is to find a minimum length Hamiltonian cycle, i.e. a minimum length cycle containing each node $v \in V$ exactly once. The 2-matching relaxation of TSP removes the subtour elimination constraints thus allowing solution to consist of several disjoint cycles. Finally, F2M problem removes integrality constraint from the 2-matching problem. For any node $v \in V$ let $\delta(v)$ denote edges originated from $v$. So LP formulation of F2M relaxation is as follows:

$$\begin{aligned} & minimize \;\; cx \\ & subject\ to: \\ & \sum_{e \in E} a_v^e x^e = 2, \; v \in V \\ & 0 \le x^e \le 1 \\ & c, x \in R^{|E|} \end{aligned} \qquad \text{(F2M)}$$

$$a_v^e = 1 \text{ if } e \in \delta(v) \text{ and } a_v^e = 0 \text{ otherwise}$$

It is well known that optimal solution $x_o$ for F2M problem is $\{0, 1/2, 1\}$ valued (Edmonds 1965). As we described in section 2, instead of solving F2M problem directly we apply gradient descent procedure to optimize $\lambda$ to make new objective function $v = c - \sum_{v \in V} \lambda_v a_v$ to be trivial:

$$\begin{aligned} & v^e < 0 \text{ for } x_o^e = 1 \\ & v^e = 0 \text{ for } x_o^e = ½ \\ & v^e > 0 \text{ for } x_o^e = 0 \end{aligned}$$

The gradient for $\lambda_v$ (for each node $v \in V$) is computed based on the length $v^e$ of adjusted edges

$$v^e = c^e - \sum_{v \in V} \lambda_v a_v^e \qquad (9)$$

that stemmed from the node $v$:

$$grad(\lambda_v) = 0.5 * ({v^e}_2 - {v^e}_3) \qquad (10)$$

where index 2 and 3 denotes the second and the third shortest edges (based on adjusted length $v^e$) originated from the node.

## 5. Computational result

5.1 Parallel implementation

The gadient (10) can be computed for each node in parallel mode as well as update of adjusted length (9) can be computed for each edge. Gradient descent algorithm described above was implemented as CUDA C/C++ program (referred further as gradient descent (GD) program ) to be executed on GPU. The input of the GD program is 2D coordinates of the nodes in TSPLIB format. As output, the optimal solution for F2M is provided: a set of edges with value 1 and a set of edges with value ½.

5.2 Benchmark

We compared the performance of our algorithm (GD) versus simplex method. For this purpose we used Soplex (Gleixner et al. 2018). SoPlex is an optimization package for solving linear programming problems (LPs) based on an advanced implementation of the primal and dual revised simplex algorithm. Each instance of F2M was converted into LP and subsequently solved by Soplex.

We would like to point out that solving F2M as LP on the extra large scale is not an optimal choice. As we already mentioned there are multiple algorithms available that can solve F2M faster (Boyd and Carr 1999). However our gradient descent procedure is applicable to a wider class of linear programs (see LP3 with condition (8)) and could be easily modified to solve more complex LP relaxations. Adaptation of the algorithms mentioned above to solve a wider class of problems is not commonly trivial.

All benchmark computations have been done on the computer with GPU Processor *GeForce RTX 2070* with 2304 cores and *AMD Processor Ryzen 5 3600 6-Core*. It is well known that the simplex method is inherently a sequential algorithm with little scope for parallelization. Soplex solver (as well as other widely used implementations of simplex method ) can use only one CPU while our GD procedure can use GPU Processor exploiting thousands of cores in parallel.

5.3 Data

We used two TSP dataset collections to benchmark our algorithm to solve fractional 2-matching problems ranging in size from several hundreds up to 200,000 nodes: VLSI data set collection (http://www.math.uwaterloo.ca/tsp/vlsi) and tsp art instances (http://www.math.uwaterloo.ca/tsp/data/art/index.html).

The table 1a provides results for some selected instances (the full results are provided in table 1B at the end ). In each case, to define a TSP graph we selected the top N closest neighbors for each node. In case of VLSI data set the N was set to be 20. In some cases this leads to suboptimal solutions (solutions that are optimal if to consider only the edges selected) as optimal F2M solutions have edges that are not included. As one can see from the table the speed up of GDP over Soplex implementation is steadily growing with the size of the F2M problem. On an extra large scale (> 20,000 nodes) the GDP performance was on the scale of a second while Soplex solving time exceeds 1000 seconds.

**Table 1a.** Benchmark results for VLSI data set instances (Selected instances). Soplex was running on one CPU (AMD Processor Ryzen 5 3600). GD was running on GeForce RTX 2070 with 2304 cores.

| **Instance ID** | **Size (Nodes)** | **Soplex (seconds)** | **GD (seconds)** | **Speed UP (ratio)** |
|---|---|---|---|---|
| dan59296 | 59296 | 2093.42 | 1.42 | 1474.2 |
| bna56769 | 56769 | 1983.17 | 4.44 | 446.7 |
| fna52057 | 52057 | 1565.3 | 1.63 | 960.3 |
| fht47608 | 47608 | 1265.53 | 0.54 | 2343.6 |
| rbz43748 | 43748 | 1077.47 | 1.35 | 798.1 |
| ics39603 | 39603 | 921.25 | 1.12 | 822.5 |
| pba38478 | 38478 | 828.31 | 2.12 | 390.7 |
| bby34656 | 34656 | 675.05 | 0.81 | 833.4 |
| ido21215 | 21215 | 255.05 | 0.25 | 1020.2 |
| xmc10150 | 10150 | 57.73 | 0.49 | 117.8 |
| dga9698 | 9698 | 55.13 | 0.16 | 344.6 |
| dlb3694 | 3694 | 8.16 | 0.11 | 74.2 |
| dcb2086 | 2086 | 2.76 | 0.09 | 30.7 |
| xit1083 | 1083 | 0.77 | 0.03 | 25.7 |
| rbu737 | 737 | 0.4 | 0.03 | 13.3 |
| pbl395 | 395 | 0.13 | 0.02 | 6.5 |
| xqf131 | 131 | 0.03 | 0.01 | 3 |

The table 2 provides results for tsp Art instances. In each case, to define a TSP graph we selected the top 10 closest neighbors for each node. In difference to VLSI data set this leads in all cases to optimal F2M solution. The trend observed for VLSI data continues for tsp Art instances.

**Table 2.** Benchmark results for TSP ART instances. Soplex was running on one CPU (AMD Processor Ryzen 5 3600). GD was running on GeForce RTX 2070 with 2304 cores.

| Instance ID | Size (Nodes) | Soplex (seconds) | GD (seconds) | Speed UP (ratio) |
|---|---|---|---|---|
| mona-liza | 100,000 | 4024.2 | 1.42 | ~3000 |
| vangogh | 120,000 | 4677.1 | 1.39 | ~3000 |
| venus | 140,000 | 5571.2 | 1.63 | ~4000 |
| pareja | 160,000 | 12659.5 | 2.34 | ~6000 |
| courbet | 180,000 | 17977.4 | 2.65 | ~7000 |
| earing | 200,000 | 26845.6 | 3.87 | ~8000 |

## 5. Discussion

Here we present a gradient descent algorithm to exploit the special structure common for many LP relaxations of combinatorial optimization problems. In difference to simplex method, gradient descent algorithm is naturally parallelised and, therefore, can exploit modern GPU. Here we implemented the gradient descent algorithm for F2M problem as CUDA C/C++ program. We demonstrated that on large scale instances gradient descent algorithm running on modern GPU achieves more than 1000 times speed up in comparison to simplex method running on a single CPU. Despite being implemented only for fractional 2-matching problems the algorithm can be easily adapted to solve more complicated relaxations (like adding subtour elimination constraint to F2M problem and so on). The program code is freely available at https://github.com/aav-antonov/F2M

**Table 1b.** Benchmark results for VLSI data set instances (Full data). Soplex was running on one CPU (AMD Processor Ryzen 5 3600). GD was running on GeForce RTX 2070 with 2304 cores.

| Instance ID | Size (Nodes) | Soplex (seconds) | GD (seconds) | Speed UP (ratio) |
|---|---|---|---|---|
| dan59296 | 59296 | 2093.42 | 1.42 | 1474.2 |
| bna56769 | 56769 | 1983.17 | 4.44 | 446.7 |
| fna52057 | 52057 | 1565.3 | 1.63 | 960.3 |
| fht47608 | 47608 | 1265.53 | 0.54 | 2343.6 |
| rbz43748 | 43748 | 1077.47 | 1.35 | 798.1 |
| ics39603 | 39603 | 921.25 | 1.12 | 822.5 |
| pba38478 | 38478 | 828.31 | 2.12 | 390.7 |
| bby34656 | 34656 | 675.05 | 0.81 | 833.4 |
| fry33203 | 33203 | 664.72 | 0.87 | 764 |
| xib32892 | 32892 | 586.82 | 0.65 | 902.8 |
| pbh30440 | 30440 | 550.02 | 0.32 | 1718.8 |
| ird29514 | 29514 | 497.17 | 1.13 | 440 |
| boa28924 | 28924 | 473.93 | 0.5 | 947.9 |
| icx28698 | 28698 | 450.35 | 0.71 | 634.3 |
| fyg28534 | 28534 | 452.16 | 0.55 | 822.1 |
| irx28268 | 28268 | 452.92 | 0.81 | 559.2 |
| bbz25234 | 25234 | 390.88 | 0.22 | 1776.7 |
| xrh24104 | 24104 | 335.57 | 0.39 | 860.4 |
| lsb22777 | 22777 | 272.84 | 0.49 | 556.8 |
| fma21553 | 21553 | 252.08 | 0.17 | 1482.8 |
| ido21215 | 21215 | 255.05 | 0.25 | 1020.2 |
| fnc19402 | 19402 | 203.11 | 0.45 | 451.4 |
| frh19289 | 19289 | 213.28 | 0.31 | 688 |
| pjh17845 | 17845 | 177.38 | 0.53 | 334.7 |
| xia16928 | 16928 | 158.96 | 0.42 | 378.5 |
| xrb14233 | 14233 | 116.64 | 0.4 | 291.6 |
| xvb13584 | 13584 | 109.82 | 0.38 | 289 |
| xmc10150 | 10150 | 57.73 | 0.49 | 117.8 |
| dga9698 | 9698 | 55.13 | 0.16 | 344.6 |
| ida8197 | 8197 | 40.34 | 0.28 | 144.1 |
| lap7454 | 7454 | 33.17 | 0.13 | 255.2 |
| bnd7168 | 7168 | 32.12 | 0.08 | 401.5 |
| xsc6880 | 6880 | 28.82 | 0.31 | 93 |
| fea5557 | 5557 | 17.82 | 0.28 | 63.6 |
| fqm5087 | 5087 | 15.03 | 0.38 | 39.6 |
| xqd4966 | 4966 | 15.02 | 0.06 | 250.3 |
| bgf4475 | 4475 | 11.73 | 0.05 | 234.6 |
| frv4410 | 4410 | 11.06 | 0.15 | 73.7 |
| bgd4396 | 4396 | 11.2 | 0.04 | 280 |
| bgb4355 | 4355 | 11.29 | 0.07 | 161.3 |
| dkf3954 | 3954 | 9.57 | 0.05 | 191.4 |
| dkc3938 | 3938 | 9.41 | 0.07 | 134.4 |
| xua3937 | 3937 | 9.22 | 0.55 | 16.8 |
| xqe3891 | 3891 | 9.07 | 0.08 | 113.4 |
| ltb3729 | 3729 | 8.5 | 0.06 | 141.7 |
| dlb3694 | 3694 | 8.16 | 0.11 | 74.2 |
| fjr3672 | 3672 | 8.55 | 0.05 | 171 |
| fjs3649 | 3649 | 8.44 | 0.19 | 44.4 |
| dhb3386 | 3386 | 6.81 | 0.09 | 75.7 |
| beg3293 | 3293 | 6.26 | 0.05 | 125.2 |
| fdp3256 | 3256 | 6.67 | 0.06 | 111.2 |
| lta3140 | 3140 | 5.75 | 0.07 | 82.1 |
| lsn3119 | 3119 | 5.87 | 0.04 | 146.8 |
| dke3097 | 3097 | 5.81 | 0.04 | 145.2 |
| pia3056 | 3056 | 5.74 | 0.25 | 23 |
| xva2993 | 2993 | 5.38 | 0.05 | 107.6 |
| dbj2924 | 2924 | 5.22 | 0.07 | 74.6 |
| lsm2854 | 2854 | 4.74 | 0.05 | 94.8 |
| irw2802 | 2802 | 4.59 | 0.04 | 114.8 |
| bch2762 | 2762 | 4.63 | 0.07 | 66.1 |
| mlt2597 | 2597 | 4.07 | 0.03 | 135.7 |
| pds2566 | 2566 | 3.87 | 0.13 | 29.8 |
| rbw2481 | 2481 | 3.69 | 0.03 | 123 |
| dea2382 | 2382 | 3.37 | 0.04 | 84.2 |
| ley2323 | 2323 | 3.03 | 0.04 | 75.8 |
| xpr2308 | 2308 | 3.32 | 0.04 | 83 |
| bck2217 | 2217 | 3.06 | 0.04 | 76.5 |
| xqc2175 | 2175 | 2.96 | 0.04 | 74 |
| bva2144 | 2144 | 2.82 | 0.14 | 20.1 |
| dcb2086 | 2086 | 2.76 | 0.09 | 30.7 |
| djb2036 | 2036 | 2.54 | 0.11 | 23.1 |
| dkd1973 | 1973 | 2.35 | 0.03 | 78.3 |
| dcc1911 | 1911 | 2.33 | 0.05 | 46.6 |
| djc1785 | 1785 | 2.07 | 0.05 | 41.4 |
| fnb1615 | 1615 | 1.58 | 0.02 | 79 |
| rby1599 | 1599 | 1.63 | 0.03 | 54.3 |
| rbv1583 | 1583 | 1.56 | 0.02 | 78 |
| fra1488 | 1488 | 1.34 | 0.03 | 44.7 |
| icw1483 | 1483 | 1.43 | 0.03 | 47.7 |
| dja1436 | 1436 | 1.32 | 0.05 | 26.4 |
| dca1389 | 1389 | 1.28 | 0.04 | 32 |
| dka1376 | 1376 | 1.25 | 0.11 | 11.4 |
| xit1083 | 1083 | 0.77 | 0.03 | 25.7 |
| pbd984 | 984 | 0.68 | 0.03 | 22.7 |
| lim963 | 963 | 0.6 | 0.02 | 30 |
| dkg813 | 813 | 0.45 | 0.02 | 22.5 |
| rbu737 | 737 | 0.4 | 0.03 | 13.3 |
| rbx711 | 711 | 0.38 | 0.05 | 7.6 |
| xql662 | 662 | 0.32 | 0.01 | 32 |
| pbm436 | 436 | 0.16 | 0.02 | 8 |
| pbn423 | 423 | 0.15 | 0.02 | 7.5 |
| pbk411 | 411 | 0.14 | 0.02 | 7 |
| pbl395 | 395 | 0.13 | 0.02 | 6.5 |
| bcl380 | 380 | 0.12 | 0.02 | 6 |
| pka379 | 379 | 0.12 | 0.03 | 4 |
| pma343 | 343 | 0.1 | 0.03 | 3.3 |
| xqg237 | 237 | 0.06 | 0.01 | 6 |
| xqf131 | 131 | 0.03 | 0.01 | 3 |